\documentclass[final,3p,10pt]{elsarticle}

\usepackage{amsmath,amssymb}
\usepackage{amsthm}

\newcommand\RR{\mathbb R}
\newcommand\PP{\mathbb P}

\newcommand\sumetilgammahtil{\sum_{\tilde{e}\in\gamma_{\tilde{h}}}}
\newcommand\TT{\mathcal T}
\newcommand\WW{\mathcal W}
\newcommand\BB{\boldsymbol B}
\newcommand\nn{{\boldsymbol n}}

\def\tm{\leavevmode\hbox{$\rm {}^{TM}$}}

\newtheorem{thm}{Theorem}
\newtheorem{lemma}[thm]{Lemma}
\newdefinition{rmk}{Remark}
\newproof{pf}{Proof}

\begin{document}

\begin{frontmatter}


 \title{A local projection stabilized method for fictitious domains}
 \author{Gabriel R. Barrenechea\corref{cor1}}
 \ead{gabriel.barrenechea@strath.ac.uk}
 \cortext[cor1]{Department of Mathematics and Statistics, University of Strathclyde, 26, Richmond Street, Glasgow G1 1XH, Scotland}
 \author{Franz Chouly\corref{cor2}}
 \ead{fchouly@univ-fcomte.fr}
 \cortext[cor2]{Laboratoire de Math\'ematiques de Besan\c con, UMR CNRS 6624, 16 route de Gray, 25030 Besan\c con}





\begin{abstract}
In this work a local projection stabilization method is proposed to solve a fictitious domain problem. The method adds a suitable fluctuation term to the formulation
thus rendering the natural space for the Lagrange multiplier stable. Stability and convergence are proved and these results are illustrated
by a numerical experiment.
\end{abstract}


\end{frontmatter}


\section{Introduction}
\label{intro}

The numerical solution of problems on smooth domains, but with complicated geometries, can be faced using different approaches, e.g., 
isoparametric elements, approximating the curved boundary by a polygonal one, etc. The problem becomes particularly complicated in the
case the domain moves or changes shape, thus forcing a constant remeshing if the shape is to be tackled in time. To avoid this, and other
complications, a fictitious domain method was proposed in \cite{GPP94} and analyzed in  \cite{GG95}. The ficticious domain approach replaces
the original problem by a mixed one on a  larger (and simpler) domain that seeks for the original variable and a Lagrange multiplier on
the physical boundary. In the analysis given in \cite{GG95}
it is proved that the combination of piecewise linear functions for the primal variable and piecewise constants for the multiplier
are inf-sup stable and convergent under the geometrical restriction that the mesh on the physical boundary 
is coarser than the mesh induced by the triangulation of the larger domain.
This is a limitation especially considering that the aforementioned intersection is needed to assemble the matrix associated to the discrete problem. 
Since then, some attempts have been made to overcome this restriction, such as 
cut elements (cf. \cite{BH10,BH11}), or XFEM approaches (cf. \cite{MBT06,HR09}).

In this work we propose a simple solution to this problem by means of a LPS-like stabilized finite element method. The starting point is the observation that the mesh for the
larger domain induces a partition of the physical boundary. 
The Lagrange multiplier space built from this partition contains a subspace such that the combination is inf-sup stable. Then, the present approach adds a fluctuation term to the 
formulation penalizing the distance between this natural finite element space and the underlying stable pair. The analysis of the new method follows then an approach related to the ones treating
minimal stabilization frameworks, such as \cite{BreFor01} and \cite{Bur08}.

The rest of the paper is organized as follows. Section \ref{prelim} introduces the notations and the problem of interest. Then, the method is presented
and its stability is proved in Section \ref{method}. Section \ref{sect-error} contains the error analysis which is illustrated by means of  a numerical experiment  in
Section \ref{numer}.

\section{Notations}\label{prelim}

We consider $\omega\subseteq \RR^2$ an open bounded domain with a Lipschitz continuous boundary $\gamma$ and outward normal vector $\nn$. To avoid technical difficulties
we will suppose that $\gamma$ is polygonal and then it's the union of $N$ straight segments $\gamma_1,\gamma_2,...,\gamma_N$. The analysis can nevertheless
be extended with minor modifications to the general case.
For $D\subset\RR^2$, the inner product on $L^2(D)$ (or $L^2(D)^2$) will be denoted by $(\cdot,\cdot)_D$.
We adopt the usual notations for Sobolev spaces. In particular, $H^{\frac{1}{2}}(\gamma)$ will be the space of traces of functions
of $H^1(\omega)$ on $\gamma$, with dual $H^{-\frac{1}{2}}(\gamma)$. The duality product on $H^{-\frac{1}{2}}(\gamma)\times H^{\frac{1}{2}}(\gamma)$
will be denoted by $\langle\cdot,\cdot\rangle_{\gamma}$. Also, for $\delta \in [0,\frac{1}{2}]$ the following space will be useful in the sequel
\begin{gather*}
\Pi_{j=1}^NH^{\delta}(\gamma_j) :=\{\xi\in L^2(\gamma) : \xi|_{\gamma_j}\in H^{\delta}(\gamma_j)\}\,.
\end{gather*}

The problem of interest reads as follows:
\begin{equation}\label{strong}
-\Delta u = \tilde{f}\qquad{\rm in}\;\omega\quad,\quad u=g\quad{\rm on}\;\gamma\,,
\end{equation}
where $\tilde{f}\in L^2(\omega)$ and $g\in H^{\frac{1}{2}}(\gamma)$.
The fictitious domain approach relies on the introduction of a larger (and simpler) domain $\Omega\supset\overline{\omega}$, an extension $f$
of $\tilde{f}$ to $\Omega$, and the solution of the following mixed problem:
{\it Find $(u,\lambda)\in \WW:=H^1_0(\Omega)\times H^{-\frac{1}{2}}(\gamma)$ such that}
\begin{gather}
(\nabla u,\nabla v)_\Omega-\langle\lambda,v\rangle_\gamma+\langle \mu,u\rangle_\gamma = (f,v)_\Omega + \langle \mu,g\rangle_\gamma\qquad\forall 
(v,\mu)\in \WW\,. \label{weak}
\end{gather}
Problems \eqref{strong} and \eqref{weak} are linked by the fact that if $(u,\lambda)$ satisfies \eqref{weak}, then $u|_\omega$ satisfies \eqref{strong} and
$\lambda$ coincides with the jump of the normal derivative of $u$ on $\gamma$ (see \cite{GPP94,GG95} for details).

To solve this weak problem, we introduce $\TT_h$, a regular triangulation of $\overline\Omega$ built using triangles $K$ with diameter $h_K$, and 
$h:=\max_{K\in\TT_h}h_K$. Let $\gamma_h$ be the partition of $\gamma$ induced by $\TT_h$, this
is, the collection of edges $e$ such that their end points are the intersections of $\gamma$ with the edges of the triangulation
$\TT_h$, plus the angular points of $\gamma$. Let also $\gamma_{\tilde{h}}$ be a partition of $\gamma$, whose vertices are also vertices of $\gamma_h$, with edges $\tilde{e}$ satisfying the following
(cf. \cite{GG95}):
there exists $C>0$ such that $3h\le|\tilde{e}|\le Ch$, for all $\tilde{e}\in\gamma_{\tilde{h}}$. Using the mesh regularity of $\TT_h$ it is easy to see that  for all $\tilde{e}\in\gamma_{\tilde{h}}$,
$\mathrm{card}\{e\in\gamma_h:e\subset\tilde{e}\}\le C$, where $C>0$ is independent of $\tilde{e}$ and $h$.

Over these partitions we define the following finite element spaces:
\begin{align*}
V_h &:= \{ v_h\in C^0(\overline{\Omega})\cap H^1_0(\Omega):v_h|_K\in\PP_1(K)\;,\;\forall K\in \TT_h\}\,,\\
\Lambda_h &:= \{ q_h\in L^2(\gamma): q_h|_e\in \PP_0(e)\;,\;\forall e\in \gamma_h\}\,,\\
\Lambda_{\tilde{h}} &:= \{ q_{\tilde{h}}\in L^2(\gamma): q_{\tilde{h}}|_{\tilde{e}}\in \PP_0(\tilde{e})\;,\;\forall \tilde{e}\in \gamma_{\tilde{h}}\}\,,
\end{align*}
and $\WW_h:=V_h\times\Lambda_h$. The pair $V_h\times\Lambda_h$ is not inf-sup stable, while, thanks to the hypothesis  on $\TT_h$ and
$\gamma_{\tilde{h}}$, the pair $V_h\times\Lambda_{\tilde{h}}$ satisfies a discrete inf-sup condition (cf. \cite{GG95}).

\section{The stabilized formulation and its stability}\label{method}

To avoid the need to use the space $\Lambda_{\tilde{h}}$, in  this work we propose the following alternative discrete problem: {\it Find $(u_h,\lambda_h)\in \WW_h$ such that:}
\begin{gather}\label{the-method}
\BB[(u_h,\lambda_h),(v_h,\mu_h)] = (f,v_h)_\Omega+\langle g,\mu_h\rangle_\gamma\qquad\forall\, (v_h,\mu_h)\in\WW_h\,,
\end{gather}
where 
\begin{gather}
\BB[(u_h,\lambda_h),(v_h,\mu_h)] = (\nabla u_h,\nabla v_h)_\Omega-\langle\lambda_h,v_h\rangle_\gamma+
\langle\mu_h,u_h\rangle_\gamma 
+\sumetilgammahtil C_s|\tilde{e}|(\lambda_h-\tilde{P}\lambda_h,\mu_h-\tilde{P}\mu_h)_{\tilde{e}}\,,\label{bilinear} 
\end{gather}
$C_s>0$, and  $\tilde{P}:L^2(\tilde{e})\to \PP_0(\tilde{e})$ stands for the orthogonal projection in $L^2(\tilde{e})$, i.e., $\tilde{P}\xi|_{\tilde{e}}:=|\tilde{e}|^{-1}(\xi,1)_{\tilde{e}}$.

Before heading to stability, we state the following preliminary result.

\begin{lemma}\label{verfurth}
There exists $\beta>0$ such that, for all $\mu_h\in\Lambda_h$
\begin{gather}
\beta\|\mu_h\|_{-\frac{1}{2},\gamma}\le \sup_{v_h\in V_h}\frac{-\langle\mu_h,v_h\rangle_\gamma}{|v_h|_{1,\Omega}}
+\left\{ \sumetilgammahtil C_s|\tilde{e}|\|\mu_h-\tilde{P}\mu_h\|^2_{0,\tilde{e}}\right\}^{\frac{1}{2}}\,.
\end{gather}
\end{lemma}

\begin{pf}
Let $\mu_h\in \Lambda_h$. Then
\begin{gather}\label{seis}
\|\mu_h\|_{-\frac{1}{2},\gamma}\le \|\mu_h-\tilde{P}\mu_h\|_{-\frac{1}{2},\gamma}+\|\tilde{P}\mu_h\|_{-\frac{1}{2},\gamma}\,.
\end{gather}
Using the definition of the norm on $H^{-\frac{1}{2}}(\gamma)$, the fact that $\tilde{P}$ is the orthogonal projection, Cauchy-Schwarz's inequality and the approximation properties
of $\tilde{P}$ (cf. \cite{EG04}) it follows that
\begin{align}
\|\mu_h-\tilde{P}\mu_h\|_{-\frac{1}{2},\gamma} &= \sup_{\xi\in H^{\frac{1}{2}}(\gamma)}
\frac{\langle \mu_h-\tilde{P}\mu_h,\xi\rangle_\gamma}{\|\xi\|_{\frac{1}{2},\gamma}}\nonumber\\
&= \sup_{\xi\in H^{\frac{1}{2}}(\gamma)}
\frac{\sumetilgammahtil (\mu_h-\tilde{P}\mu_h,\xi)_{\tilde{e}}}{\|\xi\|_{\frac{1}{2},\gamma}}\nonumber\\
&= \sup_{\xi\in H^{\frac{1}{2}}(\gamma)}
\frac{\sumetilgammahtil (\mu_h-\tilde{P}\mu_h,\xi-\tilde{P}\xi)_{\tilde{e}}}{\|\xi\|_{\frac{1}{2},\gamma}}\nonumber\\
&\le \sup_{\xi\in H^{\frac{1}{2}}(\gamma)}
\frac{\sumetilgammahtil \|\mu_h-\tilde{P}\mu_h\|_{0,\tilde{e}}\|\xi-\tilde{P}\xi\|_{0,\tilde{e}}}{\|\xi\|_{\frac{1}{2},\gamma}}\nonumber\\
&\le \sup_{\xi\in H^{\frac{1}{2}}(\gamma)}
\frac{\left\{\sumetilgammahtil |\tilde{e}|\|\mu_h-\tilde{P}\mu_h\|^2_{0,\tilde{e}}\right\}^{\frac{1}{2}}\left\{ \sum_{j=1}^N\sum_{\tilde{e}\subset\gamma_j}|\tilde{e}|^{-1}\|\xi-\tilde{P}\xi\|^2_{0,\tilde{e}}\right\}^{\frac{1}{2}}}{\|\xi\|_{\frac{1}{2},\gamma}}\nonumber\\
&\le C\,\sup_{\xi\in H^{\frac{1}{2}}(\gamma)}
\frac{\left\{\sumetilgammahtil C_s|\tilde{e}|\|\mu_h-\tilde{P}\mu_h\|^2_{0,\tilde{e}}\right\}^{\frac{1}{2}}\left\{ \sum_{j=1}^N\|\xi\|^2_{\frac{1}{2},\gamma_j}\right\}^{\frac{1}{2}}}{\|\xi\|_{\frac{1}{2},\gamma}}\nonumber\\
&\le C\,\left\{\sumetilgammahtil C_s|\tilde{e}|\|\mu_h-\tilde{P}\mu_h\|^2_{0,\tilde{e}}\right\}^{\frac{1}{2}}\,.\label{prelim1}
\end{align}
To bound the second term in \eqref{seis} we start noting that using the continuous inf-sup condition (cf. \cite{GG95}) there exists $\tilde{\beta}>0$ such that 
\begin{gather}\label{ocho} \tilde{\beta}\|\tilde{P}\mu_h\|_{-\frac{1}{2},\gamma} \,\le \sup_{v\in H^1_0(\Omega)}\frac{-\langle \tilde{P}\mu_h,v\rangle_\gamma}{|v|_{1,\Omega}}\,.\end{gather}
Next, since the pair  $V_h\times\Lambda_{\tilde{h}}$ satisfies a discrete inf-sup condition (cf. \cite{GG95}) there exists a Fortin operator
$\pi_h:H^1_0(\Omega)\to V_h$, i.e., a continuous linear operator such that $\langle\mu_{\tilde{h}},v\rangle_\gamma=\langle\mu_{\tilde{h}},\pi_h(v)\rangle_\gamma$
for all $\mu_{\tilde{h}}\in\Lambda_{\tilde{h}}$. Then, using  \eqref{ocho},  the properties of  $\pi_h$, Cauchy-Schwarz's inequality, the approximation properties of $\tilde{P}$ and the trace theorem it follows that
\begin{align}
\tilde{\beta}\|\tilde{P}\mu_h\|_{-\frac{1}{2},\gamma} &\le \sup_{v\in H^1_0(\Omega)}\frac{-\langle \tilde{P}\mu_h,\pi_h(v)\rangle_\gamma}{C|\pi_h(v)|_{1,\Omega}}\nonumber\\
&\le \sup_{v\in H^1_0(\Omega)}\frac{\langle \mu_h- \tilde{P}\mu_h,\pi_h(v)\rangle_\gamma}{C|\pi_h(v)|_{1,\Omega}}
+ \sup_{v\in H^1_0(\Omega)}\frac{-\langle \mu_h,\pi_h(v)\rangle_\gamma}{C|\pi_h(v)|_{1,\Omega}}\nonumber\\
&\le \sup_{v\in H^1_0(\Omega)}\frac{\sumetilgammahtil(\mu_h-\tilde{P}\mu_h,\pi_h(v)-\tilde{P}\pi_h(v))_{\tilde{e}}}{C|\pi_h(v)|_{1,\Omega}}
+ \sup_{v\in H^1_0(\Omega)}\frac{-\langle \mu_h,\pi_h(v)\rangle_\gamma}{C|\pi_h(v)|_{1,\Omega}}\nonumber\\
&\le C\,\sup_{v\in H^1_0(\Omega)}\frac{\left\{\sumetilgammahtil C_s|\tilde{e}|\|\mu_h-\tilde{P}\mu_h\|^2_{0,\tilde{e}}\right\}^{\frac{1}{2}}\left\{\sum_{j=1}^N\|\pi_h(v)\|^2_{\frac{1}{2},\gamma_j}\right\}^{\frac{1}{2}}}{|\pi_h(v)|_{1,\Omega}}
+ C\,\sup_{v_h\in V_h}\frac{-\langle \mu_h,v_h\rangle_\gamma}{|v_h|_{1,\Omega}}\nonumber\,,
\end{align}
and the result follows. {$\mbox{}$\hfill$\square$}
\end{pf}

We next state the main stability result for \eqref{the-method}. For this, we introduce the following mesh-dependent norm on $\WW_h$:
\begin{gather}\label{norm}
\|(v_h,\mu_h)\|^2_{\WW_h}:= |v_h|_{1,\Omega}^2+\beta^2\|\mu_h\|^2_{-\frac{1}{2},\gamma}+\sumetilgammahtil C_s|\tilde{e}|\|\mu_h-\tilde{P}\mu_h\|^2_{0,\tilde{e}}\,.
\end{gather}

\begin{thm}\label{infsup}
The bilinear form $\BB$ satisfies
\begin{gather}
\sup_{(v_h,\mu_h)\in \WW_h}\frac{\BB[(u_h,\lambda_h),(v_h,\mu_h)]}{\|(v_h,\mu_h)\|_{\WW_h}}\ge \frac{1}{6}\,\|(u_h,\lambda_h)\|_{\WW_h}\,,
\end{gather}
for all $(u_h,\lambda_h)\in \WW_h$. Hence, problem \eqref{the-method} is well-posed.
\end{thm}

\begin{pf}
Let $(u_h,\lambda_h)\in \WW_h$. From the definition of $\BB$ it easily follows that
\begin{gather}\label{elipt}
\BB[(u_h,\lambda_h),(u_h,\lambda_h)] = |u_h|^2_{1,\Omega}+\sumetilgammahtil C_s|\tilde{e}|\|\lambda_h-\tilde{P}\lambda_h\|^2_{0,\tilde{e}}\,.
\end{gather}
Next, from Lemma \ref{verfurth} there exists $w_h\in V_h$ such that $|w_h|_{1,\Omega}=\beta\|\lambda_h\|_{-\frac{1}{2},\gamma}$ and 
\begin{gather*}
\beta^2\|\lambda_h\|^2_{-\frac{1}{2},\gamma}-\beta\|\lambda_h\|_{-\frac{1}{2},\gamma}
\left\{\sumetilgammahtil C_s|\tilde{e}|\|\lambda_h-\tilde{P}\lambda_h\|^2_{0,\tilde{e}}\right\}^{\frac{1}{2}} \le -\langle \lambda_h,w_h\rangle_\gamma\,.
\end{gather*}
Then, applying Cauchy-Schwarz's and Young's inequalities we obtain
\begin{align*}
&\BB[(u_h,\lambda_h),(u_h+\frac{1}{2} w_h,\lambda_h)] = |u_h|^2_{1,\Omega}+\sumetilgammahtil C_s|\tilde{e}|\|\lambda_h-\tilde{P}\lambda_h\|^2_{0,\tilde{e}} 
+\frac{1}{2}\big[ (\nabla u_h,\nabla w_h)_\Omega - \langle\lambda_h,w_h\rangle_\gamma\big]\\
&\ge \frac{1}{2} |u_h|^2_{1,\Omega}+\sumetilgammahtil C_s|\tilde{e}|\|\lambda_h-\tilde{P}\lambda_h\|^2_{0,\tilde{e}}  
-\frac{1}{8}\beta^2\|\lambda_h\|^2_{-\frac{1}{2},\gamma}+\frac{3\beta^2}{8}\|\lambda_h\|^2_{-\frac{1}{2},\gamma}-
\frac{1}{2}\sumetilgammahtil C_s|\tilde{e}|\|\lambda_h-\tilde{P}\lambda_h\|^2_{0,\tilde{e}} \\
&= \frac{1}{2} |u_h|^2_{1,\Omega}+\frac{\beta^2}{4}\|\lambda_h\|^2_{-\frac{1}{2},\gamma}+
\frac{1}{2}\,\sumetilgammahtil C_s|\tilde{e}|\|\lambda_h-\tilde{P}\lambda_h\|^2_{0,\tilde{e}}\,,
\end{align*}
and the proof is finished  noting that $\|(u_h+\frac{1}{2} w_h,\lambda_h)\|_{\WW_h}\le \frac{3}{2}\|(u_h,\lambda_h)\|_{\WW_h}$. {$\mbox{}$\hfill$\square$}
\end{pf}

\section{Error analysis}\label{sect-error}

Supposing that $\lambda\in L^2(\gamma)$ we split the error into interpolation and discrete errors as follows
$(e_u , e_{\lambda}) :=  (u - u_h , \lambda - \lambda_h) 
 = (u - \mathcal{I}_h u , \lambda - \mathcal{J}_h \lambda)+ (\mathcal{I}_h u - u_h ,  \mathcal{J}_h \lambda - \lambda_h)
 =:  ( \eta^u , \eta^\lambda ) - ( e^h_u , e^h_\lambda)$, where $\mathcal{I}_h$ stands for the Lagrange interpolation operator, and $\mathcal{J}_h\lambda\in\Lambda_h$ is 
defined by $\mathcal{J}_h\lambda|_e:=|e|^{-1}(\lambda,1)_e$. As most of LPS-like methods,  \eqref{the-method} introduces a consistency error. Using
\eqref{weak}, \eqref{the-method} and the definition of $\BB$,   the following result is readily established.

\begin{lemma}\label{consistency}
Let us suppose that $\lambda \in L^2(\gamma)$. Then, for all $(v_h,\mu_h) \in 
\WW_h$ 
\[
\BB [ (e^u , e^\lambda) , ( v_h , \mu_h ) ]
= \sumetilgammahtil C_s | \tilde{e} | ( \lambda -  \tilde{P} \lambda , \mu_h - \tilde{P}
\mu_h )_{\tilde{e}}.
\]
\end{lemma}

Note that, though solution of \eqref{strong} can be supposed in $H^2(\omega)$, we can not expect the same regularity for the solution $u$ of \eqref{weak}, which only belongs to $H^s(\Omega)$, with $\frac32-\varepsilon \leq s \leq 2$ for any $\varepsilon>0$  (see \cite{GG95}). For the Lagrange multiplier $\lambda$, it belongs to $L^2(\gamma)$ in the worst case and to 
$\Pi_{j=1}^N H^{\frac{1}{2}}(\gamma_j)$ in the best case.
The main result of this section, namely the convergence of  method \eqref{the-method}, is stated next.

\begin{thm}
Let us suppose that $u \in 
H^s(\Omega)$ ($\frac{3}{2} - \varepsilon  \leq s \leq 2$) and
that $\lambda \in \Pi_{j=1}^N H^{\delta}(\gamma_j)$ ($0 \leq \delta \leq \frac12$). 
Then there exists a constant $C>0$,
independent of $h$, such that:
\begin{equation}
\label{cv_error}
\| (e_u , e_\lambda ) \|_{\WW_h}
\leq 
C \left ( h^{s-1} |u|_{s,\Omega} + h^{\frac12+\delta} 
\left ( \sum_{j=1}^N \| \lambda \|^2_{\delta,\gamma_j} \right )^{\frac{1}{2}}
\right ).
\end{equation}
\end{thm}

\begin{pf}
The first step is to bound the discrete error. For this, let $(w_h,t_h) \in \WW_h$ such that $\| (w_h , t_h ) \|_{\WW_h} = 1$ and the maximum on
Theorem \ref{infsup} is attained. Then, using Lemma \ref{consistency}  and Cauchy-Schwarz's inequality we arrive at
\begin{align*}
&\frac{1}{6} \| (e^h_u , e^h_\lambda ) \|_{\WW_h} \leq 
\BB [ (e^h_u , e^h_\lambda) , ( w_h , t_h ) ]\\
 = \:&\: 
- \BB [ (e_u , e_\lambda) , ( w_h , t_h ) ] + \BB [ (\eta^u , \eta^\lambda) , ( w_h , t_h ) ]\\ 
 = \:&\:
- \sumetilgammahtil C_s | \tilde{e} | ( \lambda -  \tilde{P} \lambda , t_h - \tilde{P}
t_h )_{\tilde{e}}
+ (\nabla \eta^u , \nabla w_h )_{\Omega} - \langle \eta^\lambda , w_h \rangle_{\gamma}\\
&\:
+ \langle t_h , \eta^u \rangle_{\gamma}
+ \sumetilgammahtil C_s | \tilde{e} | ( \eta^\lambda -  \tilde{P} \eta^\lambda , t_h - \tilde{P}
t_h )_{\tilde{e}}\\
\leq \:&\:
| \eta^u |_{1,\Omega} | w_h |_{1,\Omega}
+
\| \eta^{\lambda} \|_{-\frac{1}{2},\gamma} \| w_h \|_{\frac{1}{2},\gamma}
+
\| t_h \|_{-\frac{1}{2},\gamma} \| \eta^u \|_{\frac{1}{2},\gamma} 
\\
\:&\: +
\sqrt{2} \left ( 
\sumetilgammahtil C_s | \tilde{e} | \left[ \| \lambda -  \tilde{P} \lambda \|^2_{0,\tilde{e}}
+ \| \eta^\lambda -  \tilde{P} \eta^\lambda \|^2_{0,\tilde{e}}\right]
\right)^{\frac{1}{2}}
\left ( \sumetilgammahtil C_s | \tilde{e} | \| t_h -  \tilde{P} t_h \|^2_{0,\tilde{e}} 
\right)^{\frac{1}{2}} .
\end{align*}
Next, the fact that $\| (w_h , t_h ) \|_{\WW_h} = 1$, the trace Theorem and Poincar\'e's inequality lead to
\begin{gather}
\label{err1}
\| (e^h_u , e^h_\lambda ) \|_{\WW_h} \leq C \left (
| \eta^u |_{1,\Omega} 
+
\| \eta^{\lambda} \|_{-\frac{1}{2},\gamma} 
+\left ( \sumetilgammahtil C_s | \tilde{e} |\left[ \| \lambda -  \tilde{P} \lambda \|^2_{0,\tilde{e}}
+ \| \eta^\lambda -  \tilde{P} \eta^\lambda \|^2_{0,\tilde{e}}
\right]\right)^{\frac{1}{2}} \right) .
\end{gather}
Using a standard interpolation result (cf. \cite{EG04}), we have 
$ | \eta^u |_{1,\Omega} \leq C h^{s-1} | u |_{s,\Omega}$.
In addition, 
the approximation properties of $\tilde{P}$, its continuity and the approximation properties of $\mathcal{J}_h$ lead to
\begin{gather}
\label{err3}
\sum_{\tilde{e}\subset\gamma_j}\| \lambda -  \tilde{P} \lambda \|^2_{0,\tilde{e}} \leq  
C h^{2\delta} \| \lambda \|^2_{\delta,\gamma_j}\quad\textrm{and} \quad
\sum_{\tilde{e}\subset\gamma_j}\| \eta^\lambda -  \tilde{P} \eta^\lambda \|^2_{0,\tilde{e}} \leq 
C h^{2\delta} \| \lambda \|^2_{\delta,\gamma_j}\,.
\end{gather}
Finally, to bound the term $\| \eta^{\lambda} \|_{-\frac{1}{2},\gamma}$, we follow steps analogous to
\eqref{prelim1} and use $|\tilde{e}|\le Ch$ to arrive at:
\begin{align}
\label{err2}
\| \lambda - \mathcal{J}_h \lambda \|_{-\frac{1}{2},\gamma} = &
\sup_{\xi\in H^{\frac{1}{2}}(\gamma)} \frac{\sum_{j=1}^N (\lambda - \mathcal{J}_h \lambda , \xi - \mathcal{J}_h \xi )_{\gamma_j}} {\| \xi \|_{\frac{1}{2},\gamma}} 
\leq  C h^{\frac12+\delta} \left( \sum_{j=1}^N \| \lambda  \|^2_{\delta,\gamma_j} \right)^{\frac{1}{2}}\,.
\end{align}
Hence, gathering \eqref{err1}, \eqref{err2} and \eqref{err3} we obtain
\begin{equation*}
\| (e^h_u , e^h_\lambda ) \|_{\WW_h} \leq
C \left[ h^{s-1 }|u|_{s,\Omega} +  
h^{\frac12+\delta} \left ( \sum_{j=1}^N \| \lambda \|^2_{\delta,\gamma_j} \right )^{\frac{1}{2}}
\right]\,.
\end{equation*}

The interpolation error
$\| (\eta^u , \eta^\lambda ) \|_{\WW_h}$
 is bounded in the same way and we get
\begin{gather*}
\,| \eta^u |_{1,\Omega} + \beta \| \eta^{\lambda} \|_{-\frac{1}{2},\gamma} +
\left ( 
\sumetilgammahtil C_s | \tilde{e} | \| \eta^{\lambda} -  \tilde{P} \eta^{\lambda} \|^2_{0,\tilde{e}}
\right)^{\frac{1}{2}}
\leq  C \left ( h^{s-1} |u|_{s,\Omega} + 
h^{\frac12+\delta} \left ( \sum_{j=1}^N \| \lambda \|^2_{\delta,\gamma_j} \right )^{\frac{1}{2}} 
\right).
\end{gather*}
The error estimate then follows from the triangular inequality.
{$\mbox{}$\hfill$\square$}
\end{pf}

\section{A numerical experiment}\label{numer}

\begin{figure}[ht]
\centering
\includegraphics[width=0.45\textwidth]{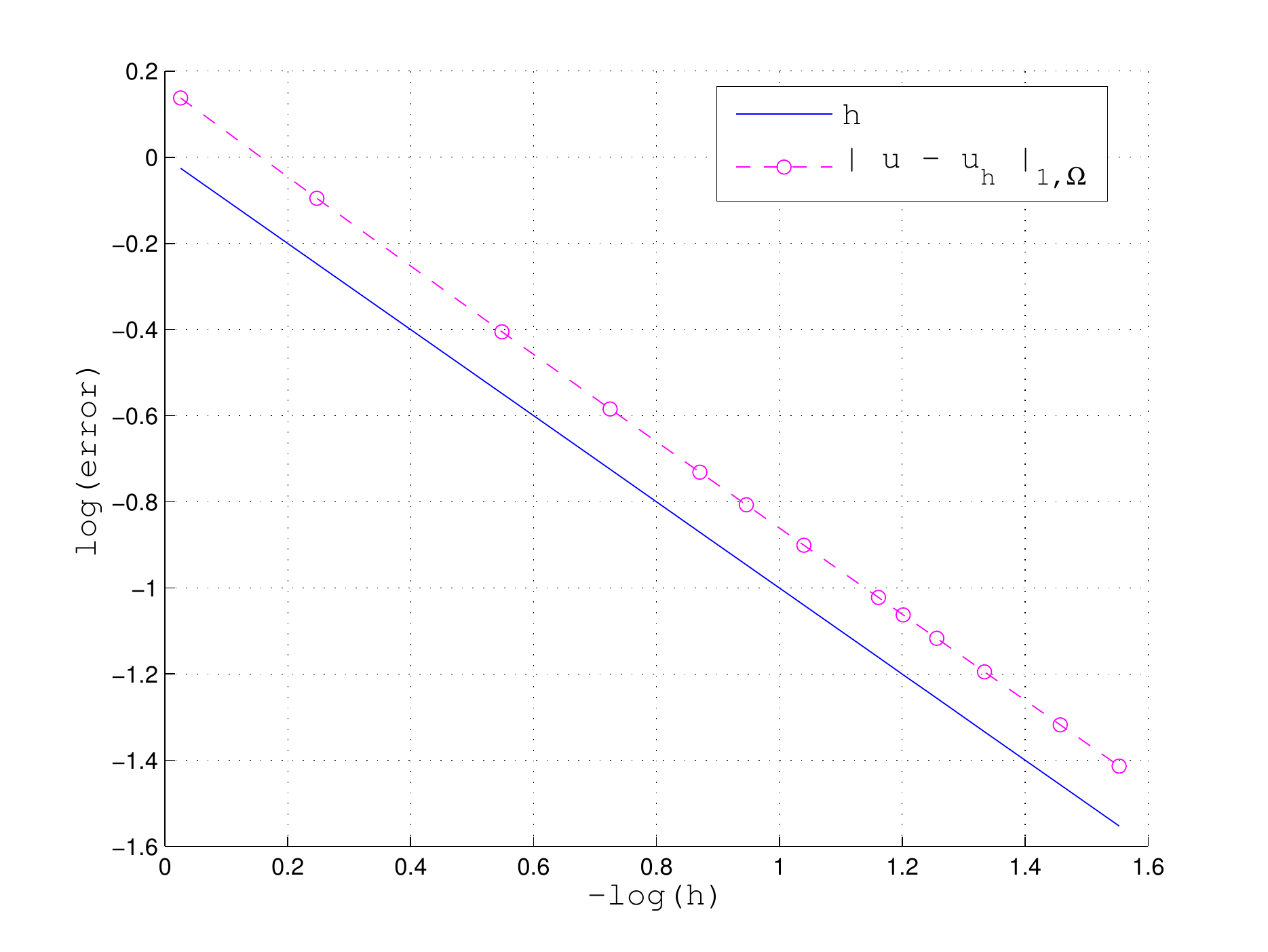} 
\includegraphics[width=0.45\textwidth]{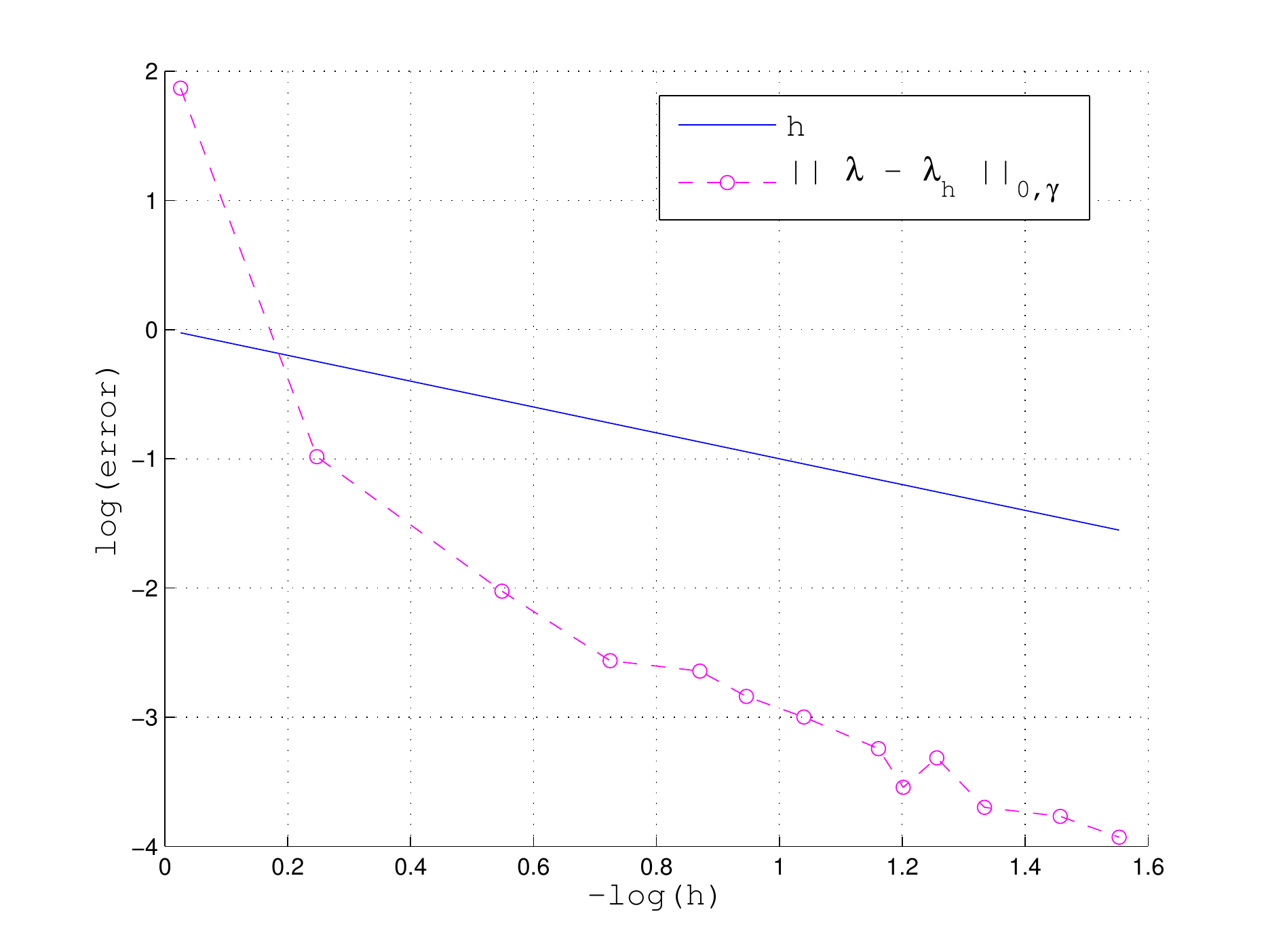}
\caption{Convergence of the method \eqref{the-method}-\eqref{bilinear}: errors $| u - u_h |_{1,\Omega}$ (left)
and $\| \lambda - \lambda_h \|_{0,\gamma}$ (right).}
\label{fig1}
\end{figure}

In order to illustrate the above theoretical results, a numerical test has been performed.
Problem \eqref{strong} has been solved using method \eqref{the-method}-\eqref{bilinear}.
We chose $\omega = [0;1]^2$, $\Omega = [-a;1+a]^2$ (with $a>0$),
and $f(x,y)=2((x+a)(1+a-x)+(y+a)(1+a-y))$ so that problem \eqref{strong} has an analytical
solution $u_a(x,y)=(x+a)(1+a-x)(y+a)(1+a-y)$. We set $g=u_a|_\gamma$ and use $a=0.5$.
A structured mesh $\TT_h$ of $\Omega$ is built, from which the boundary meshes
$\gamma_h$ and $\gamma_{\tilde{h}}$ were obtained automatically, with $\gamma_{\tilde{h}}$
satisfying $3h\le|\tilde{e}|\le 6h$. The computations have been performed
with Matlab\tm / Octave. The errors appear to be fearly independent of the value of $C_s$ in the range $0.1\le C_s\le 1000$, and then we have fixed $C_s=0.1$ in our experiments.
Convergence results are displayed Figure \ref{fig1}. Note that for the
errors $| u - u_h |_{1,\Omega}$ and $\| \lambda - \lambda_h \|_{0,\gamma}$  the optimal convergence order $O(h)$ is recovered, with a faster convergence rate for the
latter (in our case $\lambda = 0$, which helps to explain the faster convergence). This confirms
our theoretical result \eqref{cv_error}. Without stabilization ($C_s = 0$), a singular
matrix is obtained if the condition $3h\le |\tilde{e}|$ is violated. 
This confirms the necessity of the geometric condition $3h\le|\tilde{e}|\le Ch$
of \cite{GG95} without stabilization, and the interest of our stabilized formulation.

\section*{Acknowledgements}

The second author acknowledges PEPII project SIM2E and LMB (UMR CNRS 6624) for partial funding.



\bibliographystyle{model1b-num-names}
\bibliography{BC2011.bib}







\end{document}